\documentclass[leqno,english,letterpaper]{amsart}
\usepackage{verbatim}
\usepackage{graphicx}
\usepackage{amscd}
\usepackage{amsmath}
\usepackage{amsfonts}
\usepackage{amssymb}
\usepackage[all]{xy}
\theoremstyle{plain}
\newtheorem*{maintheorem}{Theorem}
\newtheorem*{intproposition}{Proposition}
\newtheorem{proposition}{Proposition}[section]
\newtheorem{theorem}[proposition]{Theorem}
\newtheorem{lemma}[proposition]{Lemma}
\newtheorem{corollary}[proposition]{Corollary}
\newtheorem{definition}[proposition]{Definition}

\theoremstyle{remark}
\newtheorem{example}[proposition]{Example}
\newtheorem{remark}[proposition]{Remark}

\numberwithin{equation}{section}

\begin{document}
\def\Z{\mathbb{Z}}
\def\ZZ{\mathbb{Z}/2\mathbb{Z}}
\def\dbar{\overline{\partial}}
\def\Dbar{\overline{D}}
\def\hom{\mathrm{Hom}}
\def\shom{\mathcal{H}\mathit{om}}
\def\ext{\mathrm{Ext}}
\def\sext{\mathcal{E}\mathit{xt}}
\def\END{\mathcal{E}\mathit{nd}}
\def\AX{\mathcal{A}_{X}}
\def\itau{\imath_{\tau}}
\def\intau{\imath_{\nabla(\tau)}}
\def\irtau{\imath_{R(\tau)}}
\def\E{\mathcal{E}}
\def\stensor{\widehat{\otimes}}
\def\CC{\mathbb{C}}
\def\O{\mathcal{O}}
\def\I{\mathcal{I}}
\def\A{\mathcal{A}}
\def\bwedge{\mbox{$\bigwedge$}}
\def\rqi{\stackrel{\sim}{\to}}
\def\lqi{\stackrel{\sim}{\gets}}
\def\cone{\mathrm{Cone}}
\def\tr{\mathrm{tr}}
\def\Tr{\mathrm{Tr}}
\def\ch{\mathrm{ch}}
\newcommand{\id}[1]{\ensuremath{\operatorname{id}_{#1}}}
\newcommand{\aff}[1]{\ensuremath{\mathbb{A}^{#1}}}
\newcommand{\proj}[1]{\ensuremath{\mathbb{P}^{#1}}}
\newcommand{\rank}{\mathrm{rank}}
\title{An Explicit Proof of the Generalized Gauss-Bonnet Formula}


\author{Henri Gillet}\email{henri@math.uic.edu}
\address{Department of MSCS, University of Illinois at Chicago, Chicago, Illinois 60607}

\thanks{The first author was supported in part by NSF Grants DMS 0100587 and DMS 0500762.}

\author{Fatih M. \"Unl\"u}
\email{unlu@imsaindy.org}
\address{Indiana Math and Science Academy, Indianapolis, Indiana 42654}



\hyphenation{grot-hen-dieck}


\begin{abstract}
In this paper we construct an explicit representative for the Grothendieck
fundamental class $ [Z] \in \ext^{r}(\O_{Z},\Omega^r_{X})$ of a complex
submanifold $Z$ of a complex manifold $X$ such that $Z$ is the zero locus of a
real analytic section of a holomorphic vector bundle $E$ of rank $r$ on $X$. To
this data we associate a super-connection $A$ on $\bigwedge^* E^{\vee}$, which
gives a ``twisted resolution'' $T^*$ of $\O_{Z}$ such that the ``generalized
super-trace'' of $\frac{1}{r!}A^{2r}$, which is a map of complexes from $T^*$ to
the Dolbeault complex $\A^{r,*}_{X}$, represents $[Z]$. One may then read off
the Gauss-Bonnet formula from this map of complexes.
\end{abstract}

\maketitle
\setcounter{section}{0}

\section*{Introduction}

If $X$ is a complex manifold, and $\tau $ is a holomorphic section,
transverse to the zero section, of the dual $E^\vee $ of  a rank $r$
holomorphic vector bundle, it is well known that the fundamental
class of the locus $Z$ of zeros of $\tau $ is equal to the top Chern
class of the bundle $E^\vee $:
\[
[Z]=c_r (E^\vee )=(-1)^rc_r (E)
\]
For Hodge cohomology, this is the fact that the image of the
Grothendieck fundamental class
$$[Z] \in \ext^{r}(\O_{Z},\Omega^r_{X})$$
under the map
$$\ext^{r}(\O_{Z},\Omega^r_{X})\to \ext^{r}(\O_{X},\Omega^r_{X})=\mathrm{H}^r(X,\Omega^r_X)$$
coincides with the top Chern class of $E^{\vee}$. Proofs of this
result tend to be indirect, i.e. they depend on the axioms for cycle
classes and Chern classes, and comparison with ``standard'' cases.

However, one may observe that the section $\tau$ gives rise to an
explicit global Koszul resolution
$$K^*(\tau)=({\bigwedge}{}^* E^{\vee},\iota_\tau)\to \mathcal{O}_Z\; ,$$
and so the theorem can be rephrased as saying that image of $[Z]$
under the map:
$$\ext^{r}(K^*(\tau),\Omega^r_{X})\to \ext^{r}(\O_{X},\Omega^r_{X})$$
induced by the isomorphism $\mathcal{O}_X\simeq K^0(\tau)$, is the
top Chern class of $E^{\vee}$. Our first result is to show that a
choice of connection $\widetilde{\nabla }$ on $E$, determines, via
Chern-Weil theory applied to superconnections, an \textit{explicit}
map of complexes from the Koszul complex $K^*(\tau)$ to the
Dolbeault complex of $\Omega_X^r$, which  represents the
Grothendieck fundamental class and the restriction of which to the
degree zero component $\mathcal{O}_X$ of the Koszul complex is
\emph{precisely} multiplication by the $r$-th Chern form of
$E^{\vee}$.

One motivation for the current paper was to obtain a better
understanding of the proof by Toledo and Tong of the
Hirzebruch-Riemann-Roch theorem in \cite{Toledo-Tong75}.  In that
paper the authors used local Koszul resolutions of the structure
sheaf of the diagonal $\Delta_X\subset X\times X$ to construct the
Grothendieck fundamental class $[\Delta_X]$, and then to compute
$\chi(X,\mathcal{O}_X)$ as the degree of the restriction of the
appropriate Kunneth-component of $[\Delta_X]$ to the diagonal. For
such a computation one needs only the existence of a ``nice''
representative of the Grothendieck fundamental class in some
neighborhood of the diagonal. However the diagonal $\Delta_X$ is
\emph{not} in general the zero set of a holomorphic section of a
vector bundle. Instead one can use the ``holomorphic exponential
map'' (see the article \cite{Kapranov99} for an exposition) to
construct, in a neighborhood of the diagonal, a \emph{real analytic}
section of $p^*(T_X)$, which vanishes exactly on the diagonal. (Here
$p:X\times X\to X$ is the projection onto the first factor.)  Thus
we are led to consider what happens if we ask only that $\tau $ be
real analytic rather than holomorphic.  In our second main result,
we use the theory of superconnections and twisted complexes in the
style of Brown \cite{Brown59}, and of Toledo and Tong (\emph{op.
cit.})
 to construct a map from the Dolbeault
resolution of $K^*(\tau)$ to that of $\Omega_X^r$ representing the
Grothendieck fundamental class and which restricts to the $r$-th
Chern form of $E^{\vee}$. An important tool in this construction is
a non-commutative version of the supertrace for endomorphisms of
Grassman algebras.

We should also remark that instead of working in the real analytic
category, one can make a very similar argument in the algebraic
category, using formal schemes.

Let us now give a more detailed outline of the paper. Recall that
the section $\tau $ gives rise to a natural Koszul resolution
$K(\tau )^{*} \to {\mathcal O}_Z $, in which $K(\tau
)^{-j}=\bigwedge ^j{\mathcal E}$. Here ${\mathcal E}$ is the sheaf
of holomorphic sections of $E$. Choose a connection $\nabla : \AX
\otimes {\mathcal E}\to {\mathcal A}^{1,0}_{X}\otimes {\mathcal E}$
of type $(1,0)$ (${\mathcal A}^{1,0}_{X}$ being the sheaf of real
analytic $(1,0)$-forms on $X)$ on ${\mathcal E}$, such that $\nabla
^2=0$. Let $\widetilde{\nabla }=\nabla +\overline
\partial $ be the associated connection. We view $\nabla $ as
acting not only on ${\mathcal E}$, but on all tensor constructions
on ${\mathcal E}$.  Then our first result is:

\begin{maintheorem}[{\bf A}]
The connection $\nabla $ and the section $\tau $ determine a map
of complexes, from the Koszul resolution $K(\tau)^*$ of ${\mathcal O}_Z $, to
the Dolbeault resolution ${\mathcal A}^{r,*}_{X}[r]$ of
$\Omega^r_{X}[r]$
\[
\psi :K (\tau )^{*} \to {\mathcal A}_{X}^{r,*}[r]
\]
the degree $-r$ component of which is $ {1 \over r!} (\intau)^{r}$, and the
degree 0 component $K(\tau)^{0}={\mathcal O}_X \to
{\A}^{r,*}_{X}[r]^{0}={\A}^{r,r}_{X}$ of which is represented by the $r$-th
Chern form of $(E^\vee ,\widetilde{\nabla })$. In general $\psi $ is given by a
linear algebra construction involving $\nabla $ and the curvature $R=[\nabla
,\overline \partial ]_{s}$ of $\widetilde{\nabla }$, and we have:
\begin{itemize}
\item The class in ${\ext}^{r}_{{\mathcal O}_{X}} ({\mathcal
O}_{Z},\Omega^r_{X})$ represented by $\psi $ is the Grothendieck fundamental
class $[Z]$.
\item The image of $[Z]$ in
${\ext}_{{\mathcal O}_X}^r ({\mathcal O}_{X},\Omega^r_{X})\simeq
{H}^{{r,r}}(X,{\mathbb C})$, is represented by the degree
zero component of $\psi $, which is equal to the $r$-th Chern form
$c_r (E^\vee,\widetilde{\nabla} )$
\end{itemize}
It follows immediately that the image of $[Z]$ in ${ H}^{r,r}(X,{\mathbb
C})$ is equal to $c_r (E^\vee )$.
\end{maintheorem}

The proof of Theorem A is contained in Section
\ref{section_Koszul_factorization}. (cf. Theorem~\ref{la13} and
Corollary~\ref{la14}).

In the second half of the paper, we extend Theorem A to the case
where $Z$ is the zero locus of a real analytic section of $E^\vee$.
It is no longer the case that $\tau $ determines a Koszul resolution
of ${\mathcal O}_Z $, but instead we get a resolution of ${\mathcal
A}^{0,*}_{X} \otimes {\mathcal O}_Z $. In order to get a complex
that is quasi-isomorphic to ${\mathcal O}_Z $, we construct a
resolution of the Dolbeault resolution ${\A}^{0,*}_{X} \otimes
{\mathcal O}_Z $ of ${\mathcal O}_Z $, by constructing a
\textit{twisted differential}, $\delta$, in the sense of Toledo and
Tong \cite{Toledo-Tong76}, on ${\mathcal A}^{0,*}_{X} \otimes
\bigwedge^{*} {\mathcal E}$.

A key tool in extending Theorem A to this situation is the notion of the
``generalized supertrace'' of an endomorphism of the exterior algebra of a
finitely generated projective module. Suppose that $V$ is a (locally) free
module of finite rank $r$ over a commutative ring $k$. Then the generalized
supertrace is a map
\[
\Tr_{\Lambda} :\mbox{End}_k (\mbox{$\bigwedge $}^{*} V)\to
\mbox{$\bigwedge $}^{*} V^\vee
\]
(c.f. Definition~\ref{la3}). If $A $ is a graded-commutative algebra over
$k$, we can extend this to a map
\[
\Tr_{\Lambda} :\mbox{End}_{A } (A \widehat{\otimes
}\mbox{$\bigwedge $} ^{*} V)\to A \widehat{\otimes
}\mbox{$\bigwedge $} ^{*} V^\vee
\]
Here $\widehat{\otimes }$ denotes the ``super'' or graded tensor
product. If $\varphi \in \mathrm{End}_{A}(\bigwedge^* V)$, then the
degree 0 component of $\Tr_{\Lambda} (\varphi) $ is the usual
super-trace of $\varphi$. The key property of $\Tr_{\Lambda} $
(which is proved in Section 3.) is:
\begin{intproposition}
Assume that $\varphi \in \mathrm{End}_{A}(\bigwedge^* V) $, and let $\delta \in
\mathrm{End}_{A}(\bigwedge^* V)$ be an $A$-linear superderivation. Then
$$
\Tr_{\Lambda}{[}\delta, \varphi {]}_{s} = {[}\delta, \Tr_{\Lambda}(\varphi)
{]}_{s} $$
\end{intproposition}

\begin{maintheorem}[{\bf B}]
Let $Z$ be a complex submanifold of $X$ such that there exists a holomorphic
vector bundle $\pi : E \rightarrow X$ and $\tau \in \Gamma(X,\mathcal{A}_X
\otimes \mathcal{E}^{\vee})$ such that $\itau: \AX \otimes \E \to \AX \otimes
\mathcal{I}_{Z} $  is surjective. Then

\begin{itemize}
\item There is a superconnection $\delta $ of type (0,1), on the
super-bundle $\bigwedge ^{*} E$, such that:
    \begin{enumerate}
    \item $\delta ^2=0,$ so $\delta $ defines a differential on
    ${\mathcal A}^{0,*}_{X} \otimes  \bigwedge ^{*} {\mathcal E}$,
    \item the
component of $\delta $ of degree $-1$ with respect to the grading on $\bigwedge
^{*} E$ is the Koszul differential $\itau $,
\item If we write $\delta
$ for the induced differential on ${\mathcal A}^{0,* }_{X}\otimes \bigwedge ^{*}
{\mathcal E}$, then the map $\bigwedge ^0{\mathcal E}={\mathcal O}_X \to
{\mathcal O}_Z $ induces a quasi-isomorphism     of complexes:
\end{enumerate} \end{itemize} \[
({\mathcal A}^{0,* }_{X} \otimes
\mbox{$\bigwedge $} ^{*} {\mathcal E},\delta ) \rqi ({\mathcal
A}^{0,* }_{X} \otimes {\mathcal O}_Z ,\overline \partial
)\;\lqi \;{\mathcal O}_Z
\]
\begin{itemize}
\item Let $R_A $ be the curvature of the superconnection
$A=\nabla +\delta $ on $\bigwedge ^{*} E$. Then the
generalized supertrace of $\frac{1}{r!}R^r_{A}$ defines a map of complexes
\[
{\mathcal A}^{0,* }_{X}\otimes
\mbox{$\bigwedge $}^{*} {\mathcal E}\to {\mathcal
A}^{r,* }_{X}[r],
\]
which, via the quasi-isomorphisms in part 1), represents the
Grothendieck fundamental class $[Z]$,
\item The image of $[Z]$ in
$H^{r,r}(X,{\mathbb C})$ is represented by the degree 0 component of
the generalized supertrace of $\frac{1}{r!}R^r_{A}$, $i.e.,$ by the super-trace of
$\frac{1}{r!}R^r_{A}$, which by Quillen \cite{Quillen85} is an $(r,r)$-form
representing the Chern character $\ch_r (\bigwedge ^{*}  E)$.
\end{itemize}
\end{maintheorem}
The proof of Theorem B is contained in Proposition~\ref{la10},
Theorem~\ref{la19}, and Corollary~\ref{la21}.
\section{Superobjects}
Throughout this paper we will use the language of
$\mathit{super}$-$\mathit{objects}$. We include here basic definitions and
properties for the convenience of the reader and to fix notation. We omit the
details and proofs, which may be found in \cite{Quillen85} and
\cite{Berline-Getzler-Vergne92}.

Let $k$ be a commutative ring with unity .
\begin{definition}
A $k$-module $V$ with a $\ZZ$-grading is called a $k$-supermodule.
\end{definition}
\begin{remark}
In the same spirit, a $\ZZ$-graded object in an additive category is called a
$\mathit{superobject}$. As realizations of this general definition, we will be
dealing with super algebras, super vector bundles on a smooth manifold, and
sheaves of superalgebras on a topological space, etc.
\end{remark}
We will write $V^{+}$ and $V^{-}$ for the degree $0 \pmod{2}$ and degree $1
\pmod{2}$ parts of $V$ and we will call them the even and the odd parts of
$V$ respectively. Let $\nu \in V$ be a homogeneous element. We say $|\nu|=0$
if $\nu \in V^{+}$ and $|\nu|=1$ if $\nu \in V^{-}$.

$\mathrm{End}_{k}(V)$ is also a $k$-supermodule with the grading
\begin{eqnarray*}
\mathrm{End}_{k}(V)^{+} &=& \hom_{k}(V^{+},V^{+}) \oplus \hom_{k}(V^{-},V^{-})
\\
\mathrm{End}_{k}(V)^{-} &=& \hom_{k}(V^{+},V^{-}) \oplus \hom_{k}(V^{-},V^{+})
\\
\end{eqnarray*}
Moreover, the algebra of endomorphisms $\mathrm{End}_{k}(V)$ is a
$k$-superalgebra with this grading. If no confusion is
likely to arise, we will suppress the mention of the ring $k$ from now on.
\begin{definition}
Let $A$ be a superalgebra. The supercommutator of
two elements of $A$ is $$
{[}a,b{]}_{s} = a b -(-1)^{|a| |b|} b a
$$
where $a$ and $b$ are homogeneous. The supercommutator is extended
bilinearly to non-homogeneous $a$ and $b$.
\end{definition}
If the supercommutator ${[}\,,{]}_{s}: A \otimes A \rightarrow A$ is the zero
map, then $A$ is called a commutative superalgebra. The exterior algebra of a
free module $M$ with the $\ZZ$-grading $\mbox{$\bigwedge$}^{+} M = \bigoplus_{{p
\ even}} \mbox{$\bigwedge$}^{p}M$ and $\mbox{$\bigwedge$}^{-} M = \bigoplus_{{p
\ odd}} \mbox{$\bigwedge$}^{p}M$ is a commutative superalgebra.

Let $V$ be finitely generated and projective. Assume that ${1 \over 2} \in
k$. Giving a $\ZZ$-grading on $V$ is equivalent to giving an involution $
\epsilon \in \mathrm{End}_{k} (V)$, that is $ \epsilon^{2} = I$. The even and
the odd parts are the eigenspaces corresponding to the eigenvalues $+1$ and $-1$
respectively. In the same fashion, the $\ZZ$-grading on $\mathrm{End}_{k}(V)$
can be given by the involution
$$
\rho(\varphi) = \epsilon \circ \varphi \circ \epsilon
$$
where $\varphi \in \mathrm{End}_{k}(V)$.
\begin{definition}
Let $\varphi \in \mathrm{End}_{k}(V)$. The supertrace of $\varphi$, denoted by
$\tr_{s}(\varphi)$, is defined to be
$$
\tr_{s}(\varphi) = \tr(\epsilon \circ \varphi)
$$
where `$\,\tr$' is the usual trace map.
\end{definition}
\begin{lemma}
The supertrace vanishes on supercommutators.
\end{lemma}
\begin{proof}
cf. \cite{Quillen85}
\end{proof}
Let $A$ and $B$ be superalgebras. We define the $\mathit{super \
tensor \ product}$ of $A$ and $B$, denoted by $A \widehat{\otimes}
B$, to be the $k$-module $A \otimes B$ with the $\ZZ$-grading
\begin{eqnarray*}
(A \widehat{\otimes} B)^{+} &=& (A^{+} \otimes B^{+}) \oplus (A^{-} \otimes
B^{-}) \\
(A \widehat{\otimes} B)^{-} &=& (A^{+} \otimes B^{-}) \oplus (A^{-} \otimes
B^{+}) \\
\end{eqnarray*}
and the algebra structure
$$
(a_{1} \otimes b_{1})(a_{2} \otimes b_{2}) = (-1)^{|b_{1}||a_{2}|} a_{1}a_{2}
\otimes b_{1}b_{2} $$
for homogeneous elements $a_{1},a_{2} \in A$ and $b_{1},b_{2} \in B$. As usual
the product is extended bilinearly.
\begin{definition}  \label{superderivation}
Let $A$ be a superalgebra and $\delta \in \mathrm{End}_{k}(A)$ be homogeneous.
We will call $\delta$ a superderivation if it satisfies the $super$-Leibniz
formula $$
\delta(a_{1}a_{2})= \delta(a_{1})a_{2} +(-1)^{|\delta||a_{1}|}a_{1}
\delta(a_{2})
$$
for homogeneous $a_{1},a_{2} \in A$. We will call a non-homogeneous element of
$\mathrm{End}_{k}(A)$ a superderivation, if its even and odd components are
superderivations.
\end{definition}
\section{Sheaves on Real Analytic Manifolds}
While we could use the $C^\infty$ Dolbeault complex in the proof of
the first main theorem, for consistency we will work with
real-analytic forms throughout this  paper.  In this section, we
shall recall the results that we need.
\begin{theorem}
Let $M$ be a real analytic manifold which is countable at infinity
and let $\mathcal{F}$ be a coherent analytic sheaf on $M$. Then
$$
H^p(M,\mathcal{F})= 0 \qquad \text{ for  } \quad  p > 0
$$
\end{theorem}
\begin{proof}
cf. Proposition 2.3 of \cite{Atiyah-Hirzebruch62}
\end{proof}
We denote the sheaf of real analytic functions on $M$ by
$\mathcal{A}_{M}$, while if $X$ is a complex manifold, we shall
write $\mathcal{A}_{X}^{p,q}$ for the sheaf of $(p,q)$-forms with
real analytic coefficients. It is a classical result (see
\cite{Dolbeault56}) that the real-analytic Dolbeault complex is a
resolution of the sheaf $\Omega_X^p$ of holomorphic $p$-forms.  It
follows from the theorem, therefore, that if $\mathcal{E}$ is a
locally free sheaf of $\mathcal{O}_X$-modules, then the cohomology
groups $H^q(X,\mathcal{E})$ may be computed as the cohomology of the
real analytic Dolbeault complex
$\mathcal{A}_{X}^{0,*}\otimes_{\mathcal{O}_X}\mathcal{E}(X)$.

\begin{corollary}  \label{la8}
Let $\mathcal{F}$ be a locally free sheaf of
$\mathcal{A}_{M}$-modules of finite rank. Then $\mathcal{F}$ is a
projective object in the category of coherent sheaves of
$\mathcal{A}_{M}$-modules.
\end{corollary}
\begin{proof}
cf. Lemma 2.7 of \cite{Atiyah-Hirzebruch62}
\end{proof}
It follows immediately that any vector bundle on a complex manifold
admits a real analytic connection, since to give such a connection
is the same as splitting the Atiyah sequence.
\begin{proposition}
Let $X_{1}$ and $X_{2}$ be complex spaces. The canonical projection
$$
\pi_{1} : X_{1} \times X_{2} \rightarrow X_{1}
$$
is flat.
\end{proposition}
\begin{proof}
cf. \cite{Fischer76} (Proposition 3.17 on page 155).
\end{proof}
\begin{corollary}
Let $X$ be a complex manifold. The sheaf $\AX$ is a flat sheaf of
$\O_{X}$-algebras.
\end{corollary}
\begin{proof}
Let $\overline{X}$ denote the complex manifold with the opposite
complex structure and $\triangle : X \rightarrow X \times
\overline{X}$ be the diagonal embedding. Let $\pi_{1} : X \times
\overline{X} \rightarrow X$ be the projection onto the first
component. Let $x \in X$ be any point. The stalks $\A_{X,x}$ and
$\O_{X \times \overline{X}, \triangle(x)}$ are canonically
isomorphic. Hence the result follows from the proposition applied to
the map $\pi_{1} : X \times \overline{X} \rightarrow X$.
\end{proof}
It follows immediately that if $\mathcal{F}$ is a coherent sheaf of
$\mathcal{O}_X$-modules, then the cohomology groups
$H^q(X,\mathcal{F})$ may be computed as the cohomology of the real
analytic Dolbeault complex
$\mathcal{A}_{X}^{0,*}\otimes_{\mathcal{O}_X}\mathcal{F}(X)$.
\section{Superconnections and the Chern Character}
Let us recall the definition and basic properties of superconnections from
\cite{Quillen85}

We assume that $X$ is a real analytic manifold of dimension $n$.
However, everything in this section applies verbatim to the smooth
case. We denote the exterior algebra of the sheaf of real analytic
differential forms on $X$, which is a sheaf of commutative
superalgebras, by $\mathcal{A}^*_{X}$. Let $E = E^{+} \oplus E^{-}$
be a real analytic super vector bundle on $X$. We will write $\E$
for the sheaf of real analytic sections of $E$, and $\A_{X}^*(\E)$
for $\A^*_{X} \widehat{\otimes}_{\AX}   \E$.
\begin{definition}   \label{superconnection}
A $\CC$-linear endomorphism $A$ of $\A_{X}^*(\E)$ of odd degree is called a
superconnection on $E$ if it satisfies the super-Leibniz rule
$$
A(\omega \otimes s)= d \omega \otimes s +(-1)^{|\omega|} \omega \wedge A(s)
$$
for local sections $\omega, s$ of $\AX^*$ and $\E$ respectively.
\end{definition}
If $X$ is an almost complex manifold and if $A$ satisfies the
following version of the super-Leibniz rule
$$
A(\omega \otimes s)= \dbar \omega \otimes s +(-1)^{|\omega|} \omega \wedge A(s)
$$
then, it is called a superconnection of type $(0,1)$ (or simply a
$(0,1)$-superconnection).

$A^{2}$ is called the curvature of the superconnection and is denoted by
$R_{A}$. The curvature of $A$ satisfies the identity
$$
R_{A}(\omega \otimes s)= \omega \wedge R_{A}(s)
$$
for local sections $\omega$ and $s$ of $\A^*_{X}$ and $\E$ respectively. Thus
$R_{A}$ can be thought as a section of the sheaf of superalgebras $\AX^*
\stensor \, \END_{\AX}(E)$ where $\END_{\AX}(E)$ denotes the sheaf of
endomorphisms of the bundle $E$.

We extend the supertrace to a map $\tr_{s} : \AX^* \stensor \,
\END_{\AX}(E) \rightarrow \AX^*$ by the formula
$$
\tr_{s}(\omega \otimes \varphi)= \omega \, \tr_{s}(\varphi)
$$
for local sections $\omega$ and $\varphi$ of $\AX^*$ and $\END_{\AX}(E)$.
\begin{proposition}
Let $n$ be a non-negative integer. The differential form $\tr_{s}(R^{n}_{A})$ is
closed and its cohomology class does not depend on the choice of the
superconnection $A$.
\end{proposition}
\begin{proof}
cf. \cite{Quillen85}
\end{proof}
\begin{theorem}   \label{la6}
The differential form
\begin{eqnarray}
\tr_{s}(\exp(R_{A}))     \label{la7}
\end{eqnarray}
represents the class $\ch(E^{+})-\ch(E^{-})$ in cohomology.
\end{theorem}
\begin{proof}
cf. \cite{Quillen85}
\end{proof}
\begin{remark}
The reader is warned that we omit the usual factor of
$({\frac{\mathrm{i}}{2\pi}})$ from \eqref{la7}, following the
convention in algebraic geometry.
 \end{remark}
\section{The Grothendieck Fundamental Class}
General references for this section are \cite{Grothendieck57} and
\cite{Altman-Kleiman70}.

Let $X$ be a compact complex manifold of dimension $n$. We denote
the sheaf of holomorphic functions and the sheaf of holomorphic $k$-forms on
$X$ by $\O_{X}$ and $\Omega_{{X}}^k$ respectively. Suppose that $\mathcal{F}$
and $\mathcal{G}$ are sheaves of $\mathcal{O}_{X}-$modules. We write
$\mathcal{H}\mathit{om}_{\mathcal{O}_{X}} \mathcal{(F,G)}$ for the sheaf of
$\mathcal{O}_{X}$-morphisms from $\mathcal{F}$ to $\mathcal{G}$ and
$\mathrm{Hom}_{\mathcal{O}_{X}}(\mathcal{F,G})$ for
$\Gamma(X,\mathcal{H}\mathit{om}_{\mathcal{O}_{X}} \mathcal{(F,G)})$. The
derived functors of $\mathcal{H}\mathit{om}_{\mathcal{O}_{X}} \mathcal{(F,G)}$
(resp. $\mathrm{Hom}_{\mathcal{O}_{X}}(\mathcal{F,G})$) will be denoted by
$\mathcal{E}\mathit{xt}^i_{\mathcal{O}_{X}}(\mathcal{F,G})$ (resp.
$\mathrm{Ext}^i_{\mathcal{O}_{X}}(\mathcal{F,G})\,)$. We simply write
$\mathcal{F}^{\vee}$ for the dual of $\mathcal{F}$.

The abelian groups $\mathrm{Ext}^i_{\mathcal{O}_{X}}(\mathcal{F,G}) $
and sheaves $\mathcal{E}\mathit{xt}^i_{\mathcal{O}_{X}}(\mathcal{F,G})$
are related by the following spectral sequence
$$
E^{i,j}_2 = H^i(X,\mathcal{E}\mathit{xt}^j_{\mathcal{O}_{X}}(\mathcal{F,G}))
\Rightarrow \mathrm{Ext}^{i+j}_{\mathcal{O}_{X}}(\mathcal{F,G}))
$$
Let $Y$ be a complex submanifold of $X$ of codimension $p$.
We denote the sheaf of ideals defining $Y$ by $\mathcal{I}$.
In this situation, one has that
\begin{eqnarray*}
\mathcal{E}\mathit{xt}^i_{\mathcal{O}_{X}}(\mathcal{O}_{Y},\Omega_{{X}}^p)&=&0
\qquad \text{ for } i < p \qquad \text{and} \\
\mathcal{E}\mathit{xt}^p_{\mathcal{O}_{X}}(\mathcal{O}_{Y},\Omega_{{X}}^p)
&=&\mbox{$\bigwedge$}^p(\mathcal{I}/\mathcal{I}^2)^{\vee} \otimes
\mathcal{O}_{Y} \otimes \Omega^p_{{X}} \\
&=&\mathcal{H}\mathit{om}_{\mathcal{O}_{X}}(\mbox{$\bigwedge$}^p
(\mathcal{I}/ \mathcal{ I}^2),\mathcal{O}_{Y} \otimes
\Omega_{{X}}^{p})
\end{eqnarray*}
All tensor products are taken over $\O_{X}$ unless stated otherwise.
It follows that the edge homomorphism
$$\mathrm{Ext}^p_{\mathcal{O}_{X}}(\mathcal{O}_{Y},\Omega_{{X}}^p)\to H^0(\mathcal{E}\mathit{xt}^p_{\mathcal{O}_{X}}(\mathcal{O}_{Y},\Omega_{{X}}^p))$$
is an isomophism, and so
$$
\mathrm{Ext}^p_{\mathcal{O}_{X}}(\mathcal{O}_{Y},\Omega_{{X}}^p)=
\mathrm{Hom}_{\mathcal{O}_{X}}(\mbox{$\bigwedge$}^p(\mathcal{I}/\mathcal{I}^2),
\mathcal{O}_{Y} \otimes \Omega_{{X}}^{p})
$$
Therefore there is a class
$[Y]$ in $ \mathrm{Ext}^p_{\mathcal{O}_{X}}(\mathcal{O}_{Y},\Omega_{{X}}^p)$
which corresponds to the homomorphism of sheaves
\begin{eqnarray*}
\mbox{$\bigwedge$}^p(\mathcal{I}/\mathcal{I}^2)\rightarrow \mathcal{O}_{Y}
\otimes \Omega^p_{{X}} \\
f_1 \wedge \cdots \wedge f_p \mapsto df_1 \wedge \cdots \wedge df_p \\
\end{eqnarray*}
The class $[Y]$ is
called the Grothendieck fundamental class of $Y$ in $X$.
\section{Koszul Factorizations}\label{section_Koszul_factorization}
In this section, we prove Theorem A of the Introduction. The proof
is contained in Proposition~\ref{la12} and Theorem~\ref{la13}.

Let $\pi:E\rightarrow X$ be a holomorphic vector bundle of rank
$r$ and let $\nabla$ be a flat real analytic connection of type
$(1,0)$ on $E$. For instance, $\nabla$ can be taken as the $(1,0)$
part of the canonical connection associated to a real analytic
hermitian structure on $E$. We write $\widetilde{\nabla}$ for
$\dbar + \nabla$. We will denote the induced connection, and the
$(1,0)$-connection on the dual bundle $E^{\vee}$ using the same
symbols. However $R$ will be used exclusively to denote the
curvature of the induced connection on $E^{\vee}$. Throughout this
section we assume that $\tau \in \Gamma(X,\mathcal{E}^{\vee})$.
Let $\itau : \AX \otimes \bigwedge^p \mathcal{E} \rightarrow \AX
\otimes \bigwedge^{p-1} \mathcal{E}$ be contraction by $\tau$ as
usual. We extend $\itau$ to an odd superderivation of the sheaf of
commutative superalgebras $\AX^* \widehat{\otimes} \bigwedge
\mathcal{E}$. Note that $\nabla(\tau) = \itau \circ \nabla +
\partial \circ \itau : \AX \otimes \mathcal{E} \rightarrow
\AX^{1,0}$ and therefore $\nabla(\tau)$ can be considered as an
element of $\Gamma(X,\AX^{1,0}\otimes \mathcal{E}^{\vee})$. We
write $\intau : \AX \otimes \bigwedge^{p} \mathcal{E} \rightarrow
\AX^{1,0} \otimes \bigwedge^{p-1} \mathcal{E}$ and $\irtau : \AX
\otimes \bigwedge^p \mathcal{E} \rightarrow \AX^{1,1} \otimes
\bigwedge^{p-1} \mathcal{E}$ for the contractions with
$\nabla(\tau) \in \Gamma (X,\AX^{1,0} \otimes \mathcal{E}^{\vee})$
and $R(\tau) \in \Gamma (X,\AX^{1,1} \otimes \mathcal{E}^{\vee})$
respectively. We extend $\intau$ (resp. $\irtau$) to an even
(resp. odd) superderivation of $\AX^* \widehat{\otimes} \bigwedge
\mathcal{E}$.

We state two facts without proof
\begin{eqnarray*}
{[}\overline{\partial},\imath_{\nabla(\tau)}{]}_s &=& \imath_{R(\tau)} \\
{[}\intau,\irtau{]}_{s} &=& 0 \\
\end{eqnarray*}
\begin{lemma}
For $1\leqslant p \leqslant r$ the following diagram is
commutative
$$
\CD \mathcal{Z}^{r-p,r-p}_{{X}} \otimes \bigwedge^p \mathcal{E}
@>{{1 \over {p!}}
\, {(\intau)}^p}>> \mathcal{A}^{r,r-p}_{{X}} \\
@V\imath_{R(\tau)}VV        @V\overline{\partial}VV \\
\mathcal{Z}^{r-p+1,r-p+1}_{{X}} \otimes \bigwedge^{p-1}
\mathcal{E} @>{{1 \over {(p-1)!}} \, {(\intau)}^{p-1}}>>
\mathcal{A}^{r,r-p+1}_{{X}} \\
\endCD
$$
where ${\mathcal{Z}}^{p,p}_{{X}}$ denote the sheaf of
$\overline{\partial}$-closed (not necessarily $\partial$-closed)
forms of type $(p,p)$.
\end{lemma}
\begin{proof}
\begin{eqnarray*}
\overline{\partial} \circ {1 \over {p!}}
\;(\imath_{\nabla(\tau)})^p &=& {1 \over {p!}} \;
{[}\overline{\partial},(\imath_{\nabla(\tau)})^p{]}_s \\ &=& {1
\over {p!}} \; \sum^{p-1}_{j=0}  (\imath_{\nabla(\tau)})^j \circ
{[}\overline{\partial},\imath_{\nabla(\tau)}{]}_s \circ
(\imath_{\nabla(\tau)})^{p-j-1} \\ &=& {1 \over {p!}} \;
\sum^{p-1}_{j=0} (\imath_{\nabla(\tau)})^j \circ \imath_{R(\tau)}
\circ (\imath_{\nabla(\tau)})^{p-j-1} \\ &=& {1 \over {p!}} \;
\sum^{p-1}_{j=0} (\imath_{\nabla(\tau)})^{p-1} \circ
\imath_{R(\tau)}  \\ &=& {1 \over {p!}} \; p \;
(\imath_{\nabla(\tau)})^{p-1} \circ \imath_{R(\tau)} = {1 \over
{(p-1)!}} \, {(\intau)}^{p-1}
\end{eqnarray*}
\end{proof}
Let $\phi_p: \bigwedge^p \mathcal{E} \rightarrow
\mathcal{H}\mathit{om}_{\mathcal{O}_{X}}(\bigwedge^{r-p}
\mathcal{E}, \bigwedge^r \mathcal{E})$ be the isomorphism given by
$$
\phi_p : \alpha \mapsto (\beta \mapsto \beta \wedge \alpha) \qquad
\text{ for} \qquad \alpha \in \mbox{$\bigwedge$}^p
\mathcal{E},\quad \beta \in \mbox{$\bigwedge$}^{r-p} \mathcal{E},
\text{ and  } 0 \leqslant p \leqslant r $$ We will identify the
sheaves  $\mathcal{H}\mathit{om}_{\mathcal{O}_{X}}(\bigwedge^{r-p}
\mathcal{E}, \bigwedge^r \mathcal{E})$ and $\bigwedge^{r-p}
\mathcal{E}^{\vee} \otimes \bigwedge^r \mathcal{E}$ via the
canonical isomorphism between them.
\begin{lemma}
The following diagram is commutative for $1 \leqslant p \leqslant
r$
$$
\CD {\bigwedge}^p \mathcal{E} @>{(\phi_p^{-1} \otimes 1) \circ
({\bigwedge}^{r-p} R \otimes 1) \circ \phi_p}>>
{\bigwedge}^p \mathcal{E} \otimes \mathcal{Z}^{r-p,r-p}_{{X}} \\
@V{\imath_{\tau}}VV  @V{\imath_{R(\tau)}}VV \\
{\bigwedge}^{p-1} \mathcal{E} @>{(\phi_{p-1}^{-1} \otimes 1) \circ
({\bigwedge}^{r-p+1} R \otimes 1) \circ \phi_{p-1}}>>
{\bigwedge}^{p-1} \mathcal{E} \otimes \mathcal{Z}^{r-p+1,r-p+1}_{{X}} \\
\endCD
$$
where $\bigwedge^p R : \bigwedge^p \mathcal{E}^{\vee} \rightarrow
\bigwedge^p \mathcal{E}^{\vee} \otimes \mathcal{Z}^{p,p}_{X}$ is
defined by $\bigwedge^p R \; (e^1 \wedge \cdots \wedge e^p)=
R(e^1) \wedge \cdots \wedge R(e^p)$ for local sections $e^1,
\ldots, e^p$ of $\mathcal{E}^{\vee}$.
\end{lemma}
\begin{proof}
The lemma is an immediate consequence of the following commutative
diagrams. (Note that $\wedge \tau$ and $\wedge R(\tau)$ denote
right multiplication by $\tau$ and $R(\tau)$ respectively).
$$
\CD \bigwedge^p \mathcal{E} @>\phi_p>> \bigwedge^{r-p}
\mathcal{E}^{\vee} \otimes \bigwedge^r \mathcal{E}
@>\bigwedge^{r-p}R \otimes 1>>  \bigwedge^{r-p} \mathcal{E}^{\vee}
\otimes \bigwedge^r \mathcal{E}
\otimes \mathcal{Z}^{r-p,r-p}_{{X}} \\
@V{\imath_{\tau}}VV   @V{\wedge \tau}VV  @V{\wedge R(\tau)}VV \\
\bigwedge^{p-1} \mathcal{E} @>\phi_{p-1}>> \bigwedge^{r-p+1}
\mathcal{E}^{\vee} \otimes \bigwedge^r \mathcal{E}
@>\bigwedge^{r-p+1}R \otimes 1>>  \bigwedge^{r-p+1}
\mathcal{E}^{\vee} \otimes \bigwedge^r \mathcal{E}
\otimes \mathcal{Z}^{r-p+1,r-p+1}_{{X}} \\
\endCD
$$
and
$$
\CD \bigwedge^{r-p} \mathcal{E}^{\vee} \otimes \bigwedge^r
\mathcal{E} \otimes \mathcal{Z}^{r-p,r-p}_{{X}}
@>{\phi_p^{-1} \otimes 1}>>   \bigwedge^p \mathcal{E} \otimes \mathcal{Z}^{r-p,r-p}_{{X}} \\
@V{\wedge R(\tau)}VV   @V{\imath_{R(\tau)}}VV \\
\bigwedge^{r-p+1} \mathcal{E}^{\vee} \otimes \bigwedge^r
\mathcal{E} \otimes \mathcal{Z}^{r-p+1,r-p+1}_{{X}}
@>{\phi_{p-1}^{-1} \otimes 1}>>   \bigwedge^{p-1} \mathcal{E}
\otimes
\mathcal{Z}^{r-p+1,r-p+1}_{{X}} \\
\endCD
$$
It is worth mentioning that $R$ can be written as a matrix of
$\dbar$-closed forms of type $(1,1)$ with respect to any given
local holomorphic framing. Consequently the image of the mapping
$\bigwedge^p R : \bigwedge^p \mathcal{E}^{\vee} \rightarrow
\bigwedge^p \mathcal{E}^{\vee} \otimes \mathcal{A}^{p,p}_{X}$ lies
in $\bigwedge^p \mathcal{E}^{\vee} \otimes \mathcal{Z}^{p,p}_{X}$.
Since $\tau$ is a holomorphic section, a similar remark applies to
the mappings $\irtau$ and $\wedge R(\tau)$.
\end{proof}
\begin{proposition} \label{la11}
The following diagram is commutative
$$
\CD \bigwedge^r \mathcal{E} @>{\imath_{\tau}}>> \bigwedge^{r-1}
\mathcal{E} @>{\imath_{\tau}}>> \cdots
@>{\imath_{\tau}}>> \mathcal{E}  @>{\imath_{\tau}}>>  \mathcal{O}_{X} \\
@V{\psi_r = {{1 \over {r!}} \, {(\intau)}^r}}VV  @V{\psi_{r-1}}VV
&& @V{\psi_1}VV  @V{\psi_0 = \det(R)}VV \\ \mathcal{A}^{r,0}_{{X}}
@>{\overline{\partial}}>>  \mathcal{A}^{r,1}_{{X}}
@>{\overline{\partial}}>> \cdots  @>{\overline{\partial}}>>
\mathcal{A}^{r,r-1}_{{X}} @>{\overline{\partial}}>>
\mathcal{A}^{r,r}_{{X}} \\ \endCD
$$
where $\psi_p ={{1 \over {p!}} \, {(\intau)}^p} \circ (\phi_p^{-1}
\otimes 1) \circ ({\bigwedge}^{r-p}R \otimes 1) \circ \phi_p $
\end{proposition}
\begin{proof}
This is an immediate result of the previous two lemmas.
\end{proof}
The symbol $f: C^{*} \rqi D^{*}$ is used to denote that $f$ is a
quasi-isomorphism.
\begin{proposition}  \label{la12}
The morphism $\psi \in \shom_{\O_{X}}(\bigwedge^{*}\E,
\A^{r,*}_{X}[r])$ represents the Grothendieck fundamental class of
$Z$ in $\sext^{r}_{\O_{X}}({\mathcal{O}_{Z}, \Omega^{r}_{X}})$.
\end{proposition}
\begin{proof}
The morphism (which is a map of complexes) $\psi \in
\shom_{\O_{X}}(K(\tau)^{*},\A^{r,*}_{X}[r])$ gives us an
$r$-cocycle, denoted by $ [ \psi ]$, in the double complex
$\shom^{*}_{\O_{X}}(K(\tau)^{.},\A^{r,.}_{X})$. We have the
quasi-isomorphisms
$$
\shom^{*}_{\O_{X}}(K(\tau)^{.},\A^{r,.}_{X}) \rqi K(\tau)^{*}
\otimes \det(\E^{\vee}) \otimes \A^{r,*}_{X} \rqi \O_{Z} \otimes
\det(\E^{\vee}) \otimes \A^{r,*}_{X}
$$
Under these maps, $ [ \psi ] $ is mapped to $ \psi_{r} :
\bigwedge^{r} \E \rightarrow \A^{r,0}_{X} \pmod{\mathcal{I}} $.

We have $\tau = \sum_{i} \alpha_{i}e^{i}$  with respect to some
local holomorphic framing $ \{ e^1, \ldots , e^r \}$ of
$E^{\vee}$. Then
\begin{eqnarray*}
\psi_{r}(e_1 \wedge \cdots \wedge e_r) &=& {1 \over r!}
\,(\intau)^{r}(e_1
\wedge \cdots \wedge e_r) \\
&=& \intau(e_1) \wedge \cdots \wedge \intau(e_r) \\
&=& \partial \alpha_{1} \wedge \cdots \wedge \partial \alpha_{r}
\pmod{\mathcal{I}}  \\
&=& d \alpha_{1} \wedge \cdots \wedge d \alpha_{r} \pmod{\mathcal{I}}  \\
\end{eqnarray*}
Then the result follows from the fact that the morphism defining
Grothendieck fundamental class of $Z$ in
$\hom_{\O_{X}}(\bigwedge^{r}(\mathcal{I} /
{\mathcal{I}}^{2}),\O_{Z} \otimes \Omega^{r}_{X})$ is mapped to
$\psi_{r}: \bigwedge^{r} \E \rightarrow \A^{r,0}_{X}
\pmod{\mathcal{I}}$ under the following sequence of
quasi-isomorphisms
$$
\bwedge^{r}(\mathcal{I} / {\mathcal{I}}^{2})^{\vee} \otimes
\Omega^{r}_{X} \rqi \O_{Z}  \otimes \det({\E^{\vee}}) \otimes
\Omega^{r}_{X} \rqi  \O_{Z} \otimes \det(\E^{\vee}) \otimes
\A^{r,*}_{X}
$$
\end{proof}
\begin{theorem}     \label{la13}
The map of complexes $\psi \in \hom_{\O_{X}}(\bigwedge^{*}\E,
\A^{r,*}_{X}[r])$ represents the Grothendieck fundamental class of
$Z$ in $\ext^{r}_{\O_{X}}({\mathcal{O}_{Z}, \Omega^{r}_{X}})$.
Moreover the image of $\psi$ in $H^{r}{(X,\Omega^{r}_{X})}$ is the
$r$-th Chern form $c_{r}(E^{\vee},\widetilde{\nabla})$.
\end{theorem}
\begin{proof}
Since $\psi$ represents the Grothendieck fundamental class
locally, it does so globally. The second result follows from the
fact that $\psi_{0} = \det(R)$ by Proposition~\ref{la11} and that
$\det(R)$ is the r-th Chern form of the pair
$(E^{\vee},\widetilde{\nabla})$.
\end{proof}
We obtain immediately
\begin{corollary}      \label{la14}
Let $\pi : E \rightarrow X$ be a holomorphic vector bundle of rank
$r$ and $\tau : X \rightarrow E^{\vee}$ be a holomorphic section
which is transverse to the zero section. If $Z$ is the complex
submanifold where $\tau$ vanishes, then the fundamental class of $Z$
in Dolbeault cohomology is represented by the $r$-th Chern form
$c_{r}(E^{\vee},\widetilde{\nabla})$.
\end{corollary}
Notice that the standard proofs of this result (for example in
\cite{Griffiths-Harris78}) implicitly use the axioms defining Chern
classes.

\section{Generalized Supertraces}\label{generalized_supertrace}
The heart of this section is Proposition~\ref{la4} which will be
used in Section 9 to construct a map of complexes that represents
the Grothendieck fundamental class.

Let $V$ be a finitely generated projective module over a
commutative ring with unity $k$, and let $V^{\vee}$ be its dual.
Let $\langle \,, \rangle : V^{\vee} \otimes_{k} V \rightarrow k$
be the pairing defined by $\langle s,t \rangle = s(t)$ for $s \in
V^{\vee}$ and $t \in V$. We extend $\langle \,, \rangle$ to a
pairing between $\bigwedge^{m} V^{\vee}$ and $\bigwedge^{m}V$ by
$$
\langle u,v \rangle = \det \langle u^{i},v_{j} \rangle
$$
where $u = u^{1} \wedge \cdots \wedge u^{m} \in \bigwedge^{m}
V^{\vee}$ and $v=v_{1} \wedge \cdots \wedge v_{m} \in
\bigwedge^{m}V$. It is easy to check that
$$
\langle u,v \rangle = ({\imath}_{u^{m}} \circ \cdots \circ
{\imath}_{u^{1}})(v_{1} \wedge \cdots \wedge v_{m} )
$$
where ${\imath}_{u^{j}}$ denotes contraction by $u^{j}$. We denote
the exterior algebra of $V$ by $\bigwedge V$ with the usual
grading for which $\bigwedge^n V$ has degree $n$. Then
$$\hom_{k}(\mbox{$\bigwedge$}V,\mbox{$\bigwedge$}V)$$
is naturally graded with
$\hom_{k}(\mbox{$\bigwedge$}^{m}V,\mbox{$\bigwedge$}^{n}V)$ having
degree $(n-m)$.
\begin{definition}   \label{la3}
Let $\varphi \in
\hom_{k}(\mbox{$\bigwedge$}V,\mbox{$\bigwedge$}V)$. If $\varphi$
has degree $ (-n)$ with $n\geqslant 0$,  we define the
$\mathit{generalized \ supertrace}$ of $\varphi$, denoted by
$\Tr_{\Lambda} (\varphi)\in \bigwedge^{n}V^{\vee}$ (as opposed to
the supertrace $\tr_{s}$), as follows:
$$
\left( \Tr_{\Lambda}(\varphi) \right) (\eta)=(-1)^{| \eta
|}\tr_{s}(\mathit{l}_{\eta}\circ \varphi)
$$
where $\mathit{l}_{\eta} \in \mathrm{End}_{k}(\bigwedge V)$ is
left multiplication by $\eta$ for some $\eta \in \bigwedge^n V$.
If $\varphi$ has positive degree, then $\Tr_{\Lambda}(\varphi)$ is
defined to be $0$.
\end{definition}
Clearly when $n=0$, we have $\Tr_{\Lambda}=\tr_{s}$.

Let $i: \bigwedge V^{\vee} \to \mathrm{End}_{k}(\bigwedge V)$ be the
inclusion defined by $i(\alpha)(\beta)=\langle \alpha, \beta
\rangle$ for $\alpha \in \bigwedge^m V^{\vee}$ and $\beta \in
\bigwedge^m V$, and $i(\alpha)(\beta)=0$ if $\beta \notin
\bigwedge^m V$ . We will often identify $\bigwedge V^{\vee}$ with
its image under $i$, and think of $\Tr_{\Lambda}(\varphi)$ as
belonging to $\mathrm{End}_{k}(\bigwedge V)$.
\begin{remark}
We have the identity
$$
\Tr_{\Lambda} \circ i = \mathrm{id}_{\bigwedge V^{\vee}}
$$
\end{remark}
\begin{remark}
Let $\tau \in V^{\vee}$ and let $\itau$ denote contraction by
$\tau$, which is a superderivation of $\bigwedge V$ of degree $-1$.
It is straightforward to check that
$$
 \Tr_{\Lambda}(\itau)=\left\{ \begin{array}{rcl}
                             0 & \text{if} & \rank \,V \geqslant 2 \\
                             \tau & \text{if} & \rank \,V=1
                             \end{array} \right.
$$
\end{remark}
An important feature of the supertrace is that it vanishes on
supercommutators. However, this is \textit{not} true of the
generalized trace $\Tr_{\Lambda}$. Instead we have:
\begin{proposition} \label{la4}
Assume that $\varphi \in \mathrm{End}_{k}(\bigwedge V) $, and let
$\delta \in \mathrm{End}_{k}(\bigwedge V)$ be a $k$-linear
superderivation. Then
$$
\Tr_{\Lambda}{[}\delta, \varphi {]}_{s} = {[}\delta,
\Tr_{\Lambda}(\varphi) {]}_{s}
$$
\end{proposition}
In order to prove the proposition we need a lemma.
\begin{lemma} \label{la5}
Assume that $\varphi \in \hom_{k}(\bigwedge V,\bigwedge V)$ is of
degree $-n \leqslant 0$, and let $\delta$ be a $k$-linear
superderivation of degree $j$  with $-n+j \leqslant 0$. Then
$$
\Tr_{\Lambda}{[}\delta,\varphi{]}_{s} = %
{[} \delta,\Tr_{\Lambda}(\varphi) {]}_{s}
$$
\end{lemma}
\begin{proof}
 Let $\eta \in
\bigwedge^{n-j}V$ be any element. Then
\begin{eqnarray*}
\Tr_{\Lambda}{[}\delta,\varphi{]}_{s} (\eta) &=& \Tr_{\Lambda}
(\delta \circ \varphi
-(-1)^{-nj}\varphi \circ \delta )(\eta) \\
&=& (-1)^{| \eta |}\tr_{s}( \mathit{l}_{\eta} \circ \delta \circ
\varphi ) - (-1)^{| \eta | -nj} \tr_{s}( \mathit{l}_{\eta} \circ
\varphi \circ \delta) \\
&& \text{by definition of}\ \Tr_{\Lambda} \\
&=& (-1)^{| \eta |}\tr_{s}( \mathit{l}_{\eta} \circ \delta \circ
\varphi ) - (-1)^{n-nj} \tr_{s}( \delta \circ \mathit{l}_{\eta}
\circ
\varphi) \\
&& \text{since}\ \tr_{s}([\delta , \mathit{l}_{\eta} \circ \varphi ]_{s})=0 \\
&=& (-1)^{n-nj+1} \tr_{s} \left( [\delta,\mathit{l}_{\eta}]_{s}
\circ \varphi \right) \\ &=& (-1)^{n-nj+1}
\tr_{s}(\mathit{l}_{\delta(\eta)}
\circ \varphi) \\
&=& -(-1)^{-nj}\Tr_{\Lambda}(\varphi)(\delta(\eta)) \quad \text{by
definition
of} \ \Tr_{\Lambda} \\
&=& [ \delta, \Tr_{\Lambda}(\varphi)]_{s}(\eta) \quad \text{since}
\ \delta \circ \Tr_{\Lambda}(\varphi) = 0 \\ \end{eqnarray*}
\end{proof}
\begin{proof}[Proof of Proposition~\ref{la4}]
We can write $\varphi = \sum_{n} \varphi_{n}$ and $\delta =
\sum_{j } \delta_{j}$ where $\varphi_{n}$ is the degree $n$
component of $\varphi$, and $\delta_{j}$ is the degree $j$
component of $\delta$. Note that a superderivation of $\bigwedge
V$ is necessarily of degree greater than or equal to $-1$.
Moreover, $\Tr_{\Lambda} {[}\delta_{j}, \varphi_{n} {]}_{s} = 0$
unless $n+j \leqslant 0$. Then
\begin{eqnarray*}
\Tr_{\Lambda}{[}\delta, \varphi {]}_{s} &=& \sum_{j \geqslant -1}
\Tr_{\Lambda}{[}\delta_{j}, \varphi {]}_{s} \\ &=& \sum_{j
\geqslant -1} \
\sum_{n \leqslant -j} \Tr_{\Lambda}{[}\delta_{j}, \varphi_{n} {]}_{s} \\
&=& \sum_{j \geqslant -1} \ \sum_{n \leqslant -j} {[}\delta_{j},
\Tr_{\Lambda}(
\varphi_{n}) {]}_{s} \quad \text{by Lemma~\ref{la5}} \\
&=& {[}\delta, \Tr_{\Lambda}(\varphi) {]}_{s} \\
\end{eqnarray*}
\end{proof}

\section{Twisted Complexes}              \label{la1}
In this section, we give a brief exposition of twisted complexes,
which were introduced by E. Brown in \cite{Brown59}. These shall be
used in Section \ref{global_resolutions} to construct ``global
resolutions" of $\O_{Z}$ by locally free $\AX$-modules. (cf.
Proposition~\ref{la10}). The reader is referred to
\cite{Toledo-Tong78} for an extensive study of the use of twisted
complexes in the duality theory of complex manifolds.

Let $\mathbf{A}$   be an abelian category.
\begin{definition}     \label{twisted_complex}
Let $M=\lbrace{M^{p,q}\rbrace}_{p,q\in\mathbb{Z} } $  be a
bigraded object in $\mathbf{A}$, and let $\delta$  be a
differential of total degree +1 on the associated graded object
$T(M)^{i}=\bigoplus_{p+q=i}M^{p,q}$. The pair $(M,\delta)$ is
called a twisted complex if $\delta$ preserves the filtration with
respect to the grading by the first degree on $M$. In this case,
$\delta$ is called the twisting differential of the pair
$(M,\delta)$.
\end{definition}
We can write $\delta= \sum_{k\geqslant0} a_k$ where $a_k \in
\prod_{p,q \in \mathbb{Z}} \hom_{\mathcal{A}}
(M^{p,q},M^{p+k,q-k+1})$. The fact that $\delta^2=0$ entails the
following, which is called the \textit{twisting cocycle
condition},
\begin{eqnarray}
\sum_{i=0}^n a_i a_{n-i} = 0  \qquad \text{for \ } n\geqslant 0.
\end{eqnarray}
Consider the cases where n=0, 1, 2
\begin{eqnarray}
a_0^2 = 0 \label{eqa} \\
a_0 a_1+a_1 a_0 = 0 \label{eqb}\\
a_0 a_2 + a_1^2 + a_2 a_0 = 0  \label{eqc}
\end{eqnarray}
$(M^{p, \ast},a_0)$ is a cochain complex for each $p \in
\mathbb{Z}$ by \eqref{eqa}. The totality of maps
$(-1)^qa_1^{p,q}:M^{p,q} \rightarrow M^{p+1,q}$ gives us a map of
complexes from $(M^{p,*},a_0)$ to $(M^{p+1,*},a_0)$ by
\eqref{eqb}. Equation~\eqref{eqc} entails that $-a_1^2$ :
$(M^{p,\ast},a_0)\rightarrow (M^{p+2,\ast},a_0)$ is chain
homotopic to the zero map, and the chain homotopy is given by $a_2
\in \prod_{p,q \in \mathbb{Z}}
\hom_{\mathcal{A}}(M^{p,q},M^{p+2,q-1})$.

By definition $(T(M)^*,\delta)$ is a filtered differential object
with the filtration given by $F^{k}= \bigoplus_{i\geqslant
k}M^{i,j}$. Hence there exists a spectral sequence with
$$
E^{p,q}_1 = H^q(M^{p,*},a_0)\Rightarrow H^{p+q}(T(M)^*,\delta)
$$
We note that $(H^*(M^{p,*},a_0),a_1)$ is a cochain complex because
of \eqref{eqc}. As a result the $E_2$ terms of the spectral
sequence are
$$
E^{p,q}_2=H^p(H^q(M^{**},a_0),a_1)
$$
\begin{remark}
The reader may observe that there is a formal similarity between
twisting cocyles and flat superconnections.
\end{remark}
\begin{example}
Every double complex (with anticommuting differentials) is a
twisted complex with $a_k=0$ for $k\geqslant 2. $
\end{example}
\begin{example}
Let $(C^{\ast},d)$ be a bounded cochain complex in $\mathbf{A}$,
and for simplicity assume that $C^k\neq 0$ only for $0\leqslant
k\leqslant n$ for some $n\in \mathbb{N}$. Suppose we are given
projective resolutions $(P^{p, \ast},\alpha_p)$ of $C^p$ for each
$p$ with augmentation maps ${\epsilon}_{p} : (P^{p,*},\alpha_{p})
\rightarrow C^p$. Then one has cochain maps $\beta_p :
(P^{p,\ast},\alpha_p)\rightarrow (P^{p+1,\ast},\alpha_{p+1})$
lifting $d:C^p\rightarrow C^{p+1}$.
$$
\CD &&       \vdots       &&         \vdots                   &&
&& \vdots          \\ &&         @V\alpha_0VV
@V\alpha_1VV           && @V\alpha_nVV     \\ 0 @>>>   P^{0,-1}
@>\beta_{-1}>>     P^{1,-1} @>\beta_{-1}>> \cdots @>\beta_{-1}>>
P^{n,-1} @>>>  0 \\ && @V\alpha_0VV @V\alpha_1VV && @V\alpha_nVV
\\ 0 @>>>   P^{0,0} @>\beta_0>> P^{1,0} @>\beta_0>> \cdots
@>\beta_0>>   P^{n,0} @>>>  0 \\ && @V\epsilon_0VV @V\epsilon_1VV
&& @V\epsilon_nVV     \\ 0 @>>>   C^0 @>d>> C^1     @>d>> \cdots
@>d>>        C^n     @>>>  0  \\ &&         @VVV @VVV &&
@VVV          \\ &&         0          && 0 &&          &&
0             \\ \endCD
$$
Now we will construct maps ${a_k}$ for $k\geqslant 0$ in order to
make $(P^{\ast\ast},\delta= \sum_{k\geqslant 0}a_k)$ a twisted
complex. (Note that $a_k=0$ for $k\geqslant n+1$).

We set $a_0 = \{ (-1)^p\alpha_p {\}}_{p \in \mathbb{Z}}$ and $a_1=
\{ \beta_q{\}}_{q \in \mathbb{Z}}$. Then the twisting cocycle
condition is satisfied for $n=0$ and $n=1$. Now assume that $a_k$
is constructed. It is easy to check that $-\sum_{i=1}^{k}
a_ia_{k-i+1}:(P^{p,\ast},a_0) \rightarrow
(P^{p+k+1,\ast}[-k+1],a_0)$ is a cochain map. If we consider
$-\sum_{i=1}^{k} a_ia_{k-i+1}$ as a map from the complex
$(P^{p,*},a_0)$ to the augmented complex $(P^{p+k+1,\ast}[-k+1]
\to C^{p+k+1}[-k+1] \to 0)$, then it is homotopic to the zero map
being a chain map from a complex of projectives to an acyclic
complex. Let $a_{k+1}$ be a chain homotopy between
$-\sum_{i=1}^{k} a_1a_{k-i+1}$ and the zero map, then one has
$$a_{k+1}a_0+a_0a_{k+1}+\sum_{i=1}^{k} a_ia_{k-i+1}=0$$ which is
equivalent to the twisting cocycle condition for $n=k+1$. This
completes the induction.

Moreover, if we endow $C^*$ with the \textit{filtration b\^ete},
i.e. the filtration defined by
$$
\sigma_{\geqslant p}(C)^{i} = \left\{ \begin{array}{rcl}
                                        0 & \text{if} & i < p \\
                                        C^{i} & \text{if} & i \geqslant p
                                        \end{array} \right.
$$
then, the augmentation map $\epsilon: T(P^{**}) \to C^*$ respects
the filtrations and is a quasi-isomorphism of the associated
graded objects. Therefore, the cochain complex
$(T(P^{\ast\ast}),\delta)$ is quasi-isomorphic to $(C^{\ast},d)$.
\end{example}
\begin{definition}     \label{la2}
Let $C$ be an object in $\mathbf{A}$. A twisted complex
$(M,\delta)$ is called a twisted resolution of $C$ if
$$
H^i(T(M)^*,\delta)= \left\{ \begin{array}{rcl}
                             C & \text{if} & i=0 \\
                             0 & \text{if} & i \neq 0
                             \end{array} \right.
$$
\end{definition}
\section{Koszul-Dolbeault Twisted Resolutions}\label{global_resolutions}
The main result of this section is Proposition~\ref{la10} wherein we
use twisted complexes to construct a global resolution of $\O_Z$ by
locally free $\AX$-modules.

Let $X$ be a compact complex manifold of dimension n. We write $\AX$ and
$\mathcal{A}^{p,q}_{{X}}$ for the sheaf of real analytic complex valued
functions and for the sheaf of real analytic differential forms of type $(p,q)$
respectively. Let $\pi: E\rightarrow X$ be a holomorphic vector bundle of rank
$r$. We will denote the sheaf of holomorphic sections of $E$ by $\mathcal{E}$.
The dual bundle and its sheaf of sections will be denoted by $E^{\vee}$ and
$\mathcal{E^{\vee}}$ respectively. \begin{definition}
A connection of type (0,1)
(or simply a (0,1)-connection) on $E$ is a $\mathbb{C}$-linear map $$
\overline{D}:\mathcal{A}_X\otimes \mathcal{E}\rightarrow \mathcal{A}_{X}^{0,1} \otimes \mathcal{E}
$$
satisfying the Leibniz rule
$$
\overline{D}(fs)=\overline{\partial}f \otimes s + f \overline{D}(s)
$$
for local sections $f$ of $\mathcal{A}_{X}$ and $s$ of $\mathcal{E}$.
\end{definition}
Let $\overline{D}$ be a $(0,1)$-connection on $E$, and let $\langle \,, \rangle$
denote the pairing between $E$ and $E^{\vee}$.
One defines a (0,1)-connection on $E^{\vee}$ (which will also be denoted by $\overline{D}$) by the formula
$$
\overline{\partial} \langle s,t \rangle =
\langle \overline{D}s,t \rangle + \langle s,\overline{D}t \rangle
$$
for local sections $s$ and $t$ of $\mathcal{A}_X \otimes \mathcal{E}$ and
$\mathcal{A}_X \otimes \mathcal{E}^{\vee}$ respectively. Finally we extend
$\Dbar$ to a $\CC$-linear superderivation of odd degree of the sheaf of
superalgebras $\AX^* \stensor \, \bigwedge \E$. As a result, $\Dbar$ is a
$(0,1)$-superconnection on the bundle $\bigwedge E$. Let $\tau$ be a global
section of $\AX \otimes \E^{\vee}$. We extend $\itau$, the contraction by
$\tau$, to an odd degree superderivation of $\AX^* \stensor \bigwedge \E$ which
acts trivially on $\A(X)$.
\begin{lemma}
Let $\tau \in \Gamma(X,\mathcal{A}_X \otimes \mathcal{E}^{\vee})$ be any
section and let $\Dbar$ be a $(0,1)$-connection on $E$ such that
$\overline{D}(\tau)=0$. The
following diagram is anti-commutative for $p,q\geqslant 0$
$$
\CD
\mathcal{A}_{X}^{0,q} \otimes {\bigwedge}^p \mathcal{E}       @>\overline{D}>>
\mathcal{A}_{X}^{0,q+1} \otimes {\bigwedge}^p \mathcal{E} \\
@V\imath_{\tau}VV    @V\imath_{\tau}VV \\
\mathcal{A}_{X}^{0,q} \otimes {\bigwedge}^{p-1} \mathcal{E}    @>\overline{D}>>
\mathcal{A}_{X}^{0,q+1} \otimes {\bigwedge}^{p-1} \mathcal{E}
\endCD
$$
\end{lemma}
\begin{proof}
$\overline{D} \circ \imath_{\tau} + \imath_{\tau} \circ \overline{D} =
[\overline{D},\imath_{\tau}]_s = \imath_{\overline{D}(\tau)}$.
\end{proof}
\begin{lemma} \label{la9}
Let $Z$ be a complex submanifold of $X$ such that there exists a
holomorphic vector bundle $\pi : E \rightarrow X$ and a section$\tau
\in \Gamma(X,\mathcal{A}_X \otimes \mathcal{E}^{\vee})$, vanishing
along $Z$, such that the induced map $\itau: \AX \otimes \E \to \AX
\otimes \mathcal{I}_{Z} $ is surjective. Then there exists a
$(0,1)$-connection $\Dbar$ on the bundle $E$ such that
$\Dbar(\tau)=0$.
\end{lemma}
\begin{proof}
We have the following diagram
$$
\CD
&&  \AX \otimes \E \\
&&  @V{\imath}_{\dbar(\tau)}VV \\
{\mathcal{A}}^{0,1}_{X} \otimes \E @>{\itau}>> {\mathcal{A}}^{0,1}_{X}
 @>p>> {\mathcal{A}^{0,1}_{X} \otimes \O_{Z}} @>>> 0 \\
\endCD
$$
There exists an $\AX$-linear map $\theta : \AX \otimes \E \rightarrow
{\mathcal{A}}^{0,1}_{X} \otimes \E $ such that $ \itau \circ \theta =
{\imath}_{\dbar(\tau)}$ since $p \circ {\imath}_{\dbar(\tau)}=0 $ and $\AX
\otimes \E$ is projective by Corollary~\ref{la8}. Then $\Dbar = \dbar -
\theta$ is the desired $(0,1)$-connection since
$$
{[}\Dbar, \itau{]}_{s} = {[}\dbar, \itau{]}_{s} - {[}\theta, \itau {]}_{s} =0
$$
\end{proof}
\begin{proposition} \label{la10}
Let $Z$ be a complex submanifold of $X$ such that there exists a holomorphic
vector bundle $\pi : E \rightarrow X$ and $\tau \in \Gamma(X,\mathcal{A}_X
\otimes \mathcal{E}^{\vee})$ such that $\itau: \AX \otimes \E \to \AX \otimes
\mathcal{I}_{Z} $  is surjective.
There is a twisted resolution (cf. Definition~\ref{la2}) $(M^{l,m},\delta =
\sum_{k \geqslant 0} a_k)$ of $\mathcal{O}_{Z}$ with
$M^{l,m}=\mathcal{A}^{0,l}_{{X}} \otimes \bigwedge^m \mathcal{E}, \quad
a_0=\imath_{\tau}, \quad a_1=\overline{D}$, and where the $a_k$ are
$\AX^*$-linear superderivations for $k \geqslant 2$.
\end{proposition}
The quasi-isomorphisms are the augmentation map $\epsilon: (T(M)^*,\delta) \to
({\A}^{0,*}_{X}\otimes \O_{Z},\dbar) $ and the inclusion $\O_{Z}
\to ({\A}^{0,*}_{X}\otimes \O_{Z},\dbar)$.
\begin{proof}
The fact that $({\itau})^{2}=0$ and that $\Dbar$ and $\itau$ anticommute implies
that the twisting cocycle condition (defined in Section~\ref{la1}) is
satisfied for $n=0$ and $n=1$. We will construct superderivations $a_{k}, \quad
k \geqslant 2$ satisfying the twisting cocycle condition by induction.
$$
\CD
\mathcal{A}_{X}  @>\overline{\partial}>>  \mathcal{A}^{0,1}_{{X}}  @>\overline{\partial}>>  \cdots
@>\overline{\partial}>>
\mathcal{A}^{0,n}_{{X}} \\
@A\imath_{\tau}AA  @A\imath_{\tau}AA  && @A\imath_{\tau}AA \\
\mathcal{A}_{X} \otimes \mathcal{E}  @>\overline{D}>>  \mathcal{A}^{0,1}_{{X}}
\otimes \mathcal{E} @>\overline{D}>>  \cdots  @>\overline{D}>>
\mathcal{A}^{0,n}_{{X}} \otimes \mathcal{E}\\
@A\imath_{\tau}AA  @A\imath_{\tau}AA  && @A\imath_{\tau}AA \\
\vdots && \vdots  &&  && \vdots \\
@A\imath_{\tau}AA  @A\imath_{\tau}AA  && @A\imath_{\tau}AA \\
\mathcal{A}_{X} \otimes \bigwedge^{r-1} \mathcal{E} @>\overline{D}>>
\mathcal{A}^{0,1}_{{X}} \otimes \bigwedge^{r-1} \mathcal{E}
@>\overline{D}>> \cdots  @>\overline{D}>>  \mathcal{A}^{0,n}_{{X}} \otimes \bigwedge^{r-1}
\mathcal{E}\\
@A\imath_{\tau}AA  @A\imath_{\tau}AA  && @A\imath_{\tau}AA \\
\mathcal{A}_{X} \otimes \bigwedge^r \mathcal{E} @>\overline{D}>>
\mathcal{A}^{0,1}_{{X}} \otimes \bigwedge^r \mathcal{E}
@>\overline{D}>> \cdots  @>\overline{D}>>  \mathcal{A}^{0,n}_{{X}} \otimes \bigwedge^r \mathcal{E}\\
\endCD
$$
We have the diagram
$$
\CD
&&  {\A}_{X} \otimes \E @>{\itau}>> \AX \\
&&  @V{-a^{2}_{1}}VV   @V{0}VV \\
{\A}^{0,2}_{X} \otimes \bigwedge^{2} \E @>{\itau}>>  {\A}^{0,2}_{X} \otimes
\E @>{\itau}>>  {\A}^{0,2}_{X} \\
\endCD
$$
There exists an $\AX$-linear map $a_{2} : {\A}_{X} \otimes \E
\rightarrow {\A}^{0,2}_{X} \otimes \bigwedge^{2} \E$ such that $-a^{2}_{1} =
\itau \circ a_{2}$ since ${\A}^{0,q}_{X} \otimes \E$ is projective by
Lemma~\ref{la8}. We extend $a_{2}$ to an odd superderivation of
$\AX^* \stensor \bigwedge \E$ which acts trivially on $\AX^*$. The twisting
cocycle condition for $n=2$ is satisfied since both $-a^{2}_{1}$ and $a_{2}
\circ \itau + \itau \circ a_{2}$ are superderivations that act trivially on
$\AX^*$ and agree on $\AX \otimes \E$.

Now assume that $a_{k}$ is constructed. Thus we have the diagram
$$
\CD
&&  {\A}_{X} \otimes \E @>{\itau}>> \AX \\
&&  @V{\mu}VV   @V{0}VV \\
{\A}^{0,k+1}_{X} \otimes \bigwedge^{k+1} \E @>{\itau}>>  {\A}^{0,k+1}_{X}
\otimes \bigwedge^{k}\E @>{\itau}>>  {\A}^{0,k+1}_{X} \otimes
\bigwedge^{k-1} \E   \\
\endCD
$$
where $\mu = - \sum_{i=1}^{k} a_{i}a_{k-i+1}$. It is straightforward to check
that $\mu$ is $\AX$-linear. Hence there exists a map $a_{k+1}
: {\A}_{X} \otimes \E \rightarrow {\A}^{0,k+1}_{X} \otimes
\bigwedge^{k+1} \E$ such that $ \mu = \itau \circ a_{k+1}$. The extension of
$a_{k+1}$ to an odd superderivation of $\AX^* \stensor \bigwedge \E$
which acts trivially on $\AX^*$ satisfies the twisting cocycle condition by
applying the argument used in the previous paragraph. Note that $a_{k}=0$ for $k
\geqslant n+1$, so the induction ends after a finite number of steps.

Let $\epsilon$ denote the augmentation map from the twisted complex
$(T^*=T(M)^{*},\delta)$ to the complex $({\A}^{0,*}_{X}\otimes \O_{Z},\dbar)$.
If we endow the latter complex with \textit{filtration b\^ete}, then $\epsilon$
becomes a map of filtered complexes which is a quasi-isomorphism of the
associated graded objects. Hence, $(T^*,\delta)$ is quasi-isomorphic to
$({\A}^{0,*}_{X}\otimes \O_{Z},\dbar)$, which in turn is quasi-isomorphic to
$\mathcal{O}_{Z}$.
\end{proof}
Note that $\delta$ is a flat superconnection of type $(0,1)$ on the superbundle
$\bigwedge E$.
\section{Koszul Factorizations II}
In this section, we will construct a map from the twisted complex
$T^*=T(M)^*$ of Proposition~\ref{la10} to the Dolbeault complex
$\A^{r,*}_{X}[r]$ by using the generalized supertraces of Section
\ref{generalized_supertrace}. The precise argument is contained in
Corollary~\ref{la16}. In the next section, we will prove that this
map represents the Grothendieck fundamental class of $Z$ in $X$.

We extend the generalized trace (cf. Definition~\ref{la3}) to a map
$$
\Tr_{\Lambda}: \AX^* \stensor \ \END_{\O_{X}}(\bwedge E) \rightarrow \AX^*
\stensor \bwedge {\E}^{\vee}
$$
by the formula
$$
\Tr_{\Lambda}(\omega \otimes \varphi) = \omega \, \Tr_{\Lambda}(\varphi)
$$
for local sections $\omega$ and $\varphi$ of $\AX^*$ and
$\END_{\O_{X}}(\bigwedge E)$ respectively.
\begin{proposition}
Let $\varphi$ be a section of the sheaf of superalgebras $ \AX^*
\stensor \ \END_{\O_{X}}(\bigwedge E) $ and $\delta$ be the twisting
differential of Proposition~\ref{la10}. Then
$$
\Tr_{\Lambda}{[} \delta, \varphi {]}_{s} = {[} \delta, \Tr_{\Lambda}(\varphi){]}_{s}
$$
\end{proposition}
\begin{proof}
We observe that, for local sections ${\omega}_{1} \otimes {\varphi}_{1},
{\omega}_{2} \otimes {\varphi}_{2}$ of $\AX^* \stensor \ \END_{\O_{X}}(\bwedge
E)$, one has that
$$
{[}{\omega}_{1} \otimes {\varphi}_{1},{\omega}_{2} \otimes
{\varphi}_{2}{]}_{s}=(-1)^{|\varphi_{1}||\omega_{1}|} {\omega_{1}} \wedge
{\omega_{2}} \otimes {[} {\varphi_{1}},{\varphi_{2}}{]}_{s}
$$
This follows from a straightforward computation and the fact that $\AX^*$ is
supercommutative. Therefore
$$
\Tr_{\Lambda}{[}{\omega}_{1} \otimes {\varphi}_{1},{\omega}_{2} \otimes
{\varphi}_{2}{]}_{s} = (-1)^{|\varphi_{1}||\omega_{1}|} {\omega_{1}} \wedge
{\omega_{2}} \otimes \Tr_{\Lambda} {[} {\varphi_{1}},{\varphi_{2}}{]}_{s}
$$
Assume that $\tilde{\delta}$ is a section of $\AX^* \stensor \
\END_{\O_{X}}(\bwedge E)$ which is a superderivation. Since the
supercommutator and $\Tr_{\Lambda}$ are additive, we may assume
(without any loss of generality) that $\tilde{\delta}  =
{\omega}_{1} \otimes {\varphi}_{1}$ and $\varphi= {\omega}_{2}
\otimes {\varphi}_{2}$. Then \begin{eqnarray*} \Tr_{\Lambda}{[}
\tilde{\delta}, \varphi{]}_{s} &=& (-1)^{|\varphi_{1}||\omega_{1}|}
{\omega_{1}} \wedge {\omega_{2}} \otimes \Tr_{\Lambda} {[}
{\varphi_{1}},{\varphi_{2}}{]}_{s} \\ &=&
(-1)^{|\varphi_{1}||\omega_{1}|} {\omega_{1}} \wedge {\omega_{2}}
\otimes  {[} {\varphi_{1}},\Tr_{\Lambda} ({\varphi_{2}}){]}_{s}
\quad  \text{by Proposition~\ref{la4}} \\ &=& {[}{\omega}_{1}
\otimes {\varphi}_{1},
{\omega_{2}} \otimes \Tr_{\Lambda}( {\varphi}_{2}) {]}_{s} \\
&=& {[} \tilde{\delta} , \Tr_{\Lambda}(\varphi){]}_{s}
\end{eqnarray*}
We further observe that $\delta - \dbar$ is a section of $\AX^* \stensor \
\END_{\O_{X}}(\bwedge E)$ which is a sum of superderivations. As a result we have
$$
\Tr_{\Lambda}{[}\delta -\dbar,\varphi{]}_{s} = {[}\delta -\dbar ,
\Tr_{\Lambda}(\varphi){]}_{s}
$$
We finally observe that
$$
{[}\dbar,\varphi{]}_{s}={[} \dbar, {\omega}_{2} \otimes {\varphi}_{2}{]}_{s}=
\dbar{\omega}_{2} \otimes {\varphi}_{2}
$$
and
$$
{[}\dbar ,\Tr_{\Lambda}(\varphi){]}_{s}={[}\dbar , \Tr_{\Lambda}({\omega}_{2} \otimes
{\varphi}_{2}){]}_{s}={[}\dbar , {\omega}_{2} \otimes \Tr_{\Lambda}(
{\varphi}_{2}){]}_{s}= \dbar {\omega}_{2} \otimes \Tr_{\Lambda}({\varphi}_{2})
$$
Then the assertion follows from the equation
$$
\Tr_{\Lambda}{[}\dbar,\varphi{]}_{s} = {[}\dbar ,
\Tr_{\Lambda}(\varphi){]}_{s}
$$
\end{proof}
Recall that the differential $\delta$ is a flat superconnection of type $(0,1)$
on the super vector bundle $\bwedge E$. If we let ${A}= \nabla + \delta$ (recall
that $\nabla$ is a flat $(1,0)$-connection on $E$),
then $A$ is a superconnection on $\bwedge E$. The curvature of $A$, denoted by
$R_A$, is given by the formula
$$
R_A = A^2 = (\nabla + \delta)^2 = \nabla^2 + \nabla \circ \delta + \delta \circ
\nabla + \delta^2 = {[}\nabla,\delta{]}_{s}
$$
since $\nabla^2 = \delta^2 =0$.
\begin{corollary} \label{la15}
Let $ \psi = \frac{1}{r!} R_{A}^{r}$. Then
$$
{[}\delta, \Tr_{\Lambda}(\psi){]}_{s} =0
$$
\end{corollary}
\begin{proof}
This follows from the fact that ${[} \delta, R_{A}{]}_{s}=0$ and the
proposition.
\end{proof}
\begin{corollary}    \label{la16}
Let $ \psi =\frac{1}{r!} R_{A}^{r}$. Then
$$
\dbar \circ \Tr_{\Lambda}(\psi) = \Tr_{\Lambda} (\psi) \circ \delta
$$
In other words, $\Tr_{\Lambda}(\psi)$ is a cochain map from the twisted complex
$(T^*=T(M)^*,\delta)$ of Proposition~\ref{la10} (which is quasi-isomorphic to
$\O_{Z}$) to the Dolbeault complex $(\A^{r,*}_{X}[r],\dbar)$ (which is
quasi-isomorphic to $\Omega^{r}_{X}[r]$).
\end{corollary}
\begin{proof}
We have
\begin{eqnarray*}
[ \delta, \Tr_{\Lambda}(\psi) {]}_{s} &=& \sum_{j \geqslant 0} \ \sum_{m-n \geqslant
j-1} [a_{j},\Tr_{\Lambda}(\psi_{m,n}){]}_{s} \\
&=& \sum_{j \geqslant 0} \ \sum_{m-n \geqslant j-1} (a_{j} \circ
\Tr_{\Lambda}(\psi_{m,n}) - \Tr_{\Lambda}(\psi_{m,n}) \circ a_{j}) \\
&& \text{ since  } \Tr_{\Lambda}(\psi_{m,n}) \text{ are of even degree} \\
&=& \dbar \circ \Tr_{\Lambda}(\psi) - \sum_{j \geqslant 0} \ \sum_{m-n \geqslant j-1}
\Tr_{\Lambda}(\psi_{m,n}) \circ a_{j}  \\
&=& \dbar \circ \Tr_{\Lambda}(\psi) - \Tr_{\Lambda}(\psi) \circ \delta \\
\end{eqnarray*}
Then the assertion follows from Corollary~\ref{la15}
\end{proof}
Corollary~\ref{la16} (combined with the Lemmas \ref{la17} and
\ref{la20}) can be seen as a generalization of Proposition~\ref{la11} to the
real-analytic case.
\section{Comparison with the Grothendieck Class}
As a result of Corollary~\ref{la16}, $\Tr_{\Lambda}(\psi)$ gives us an element in
$\hom_{\O_{X}}(T^*,\A^{r,*}_{X}[r])$, and therefore a class in
$\ext^{r}_{\O_{X}}(T^*,\A^{r,*}_{X})$. We can identify
$\ext^{r}_{\O_{X}}(T^*,\A^{r,*}_{X})$ with the group
$\ext^{r}_{\O_{X}}(\O_{Z},\Omega^{r}_{X})$ since $T^*$ and $\A^{r,*}_{X}$ are
quasi-isomorphic to $\O_{Z}$ and $\Omega^{r}_{X}$ respectively. Now we shall
prove that the class of $\Tr_{\Lambda}(\psi)$ in
$\ext^{r}_{\O_{X}}(\O_{Z},\Omega^{r}_{X})$, denoted by $ [ \Tr_{\Lambda}(\psi) ] $, is
the Grothendieck fundamental class. But we need two preliminary lemmas first.

\begin{lemma}    \label{la17}
Let $ \psi = \frac{1}{r!}R_{A}^{r}$. We write $\psi = \sum {\psi}_{m,n}$ where
${\psi}_{m,n} \in \A(X) \otimes {\hom}_{\AX}(\bigwedge^{m} \E, \bigwedge^{n}
\E)$. Then
$$
\psi_{r,0} = \frac{1}{r!}(\intau)^{r}
$$
\end{lemma}
\begin{proof}
We have $R_{A}= {[}\nabla, \delta {]}_{s}= \intau +\sum_{k \geqslant 1}
\nabla(a_{k})$. In this sum, the only summand that lowers the Koszul degree is
the term $\intau$. Therefore, the only term that lowers the Koszul degree by $r$
in $\frac{1}{r!}R_{A}^{r}$ is $\frac{1}{r!}(\intau)^{r}$.
\end{proof}
Therefore, we have the following equality for the restriction of $\Tr_{\Lambda}(\psi)$
to ${\A_{X} \otimes \bwedge^{r} \E}$
$$
\Tr_{\Lambda}(\psi){|}_{\A_{X} \otimes \bwedge^{r} \E} = \frac{1}{r!}(\intau)^{r}
$$
\begin{lemma}    \label{la18}
Let $A^*, B^*$, and $C^*$ be cochain complexes; and $ f: A^* \rightarrow B^*$
and $g: B^* \rightarrow C^*$ be maps of complexes such that the composition $g
\circ f$ is homotopic to the zero map. There exists a map $l(f): A^*
\rightarrow \cone(g)^* [-1]$ such that the following triangle is commutative
$$
\xymatrix{
A^* \ar[r]^{l(f) {\quad}} \ar[dr]_f &  \cone(g)^* [-1] \ar[d]^{\Pr} \\
   & B^* \\
}
$$
where $\Pr : \cone(g)^* [-1] \to B^*$ is the projection map.
\end{lemma}
\begin{proof}
Exercise.
\end{proof}
\begin{theorem}      \label{la19}
 $ [ \Tr_{\Lambda}(\psi) ] \in \ext^{r}_{\O_{X}}(\O_{Z},\Omega^{r}_{X})$ is the
Grothendieck fundamental class.
\end{theorem}
\begin{proof}
Since $\ext^{r}_{\O_{X}}(\O_{Z},\Omega^{r}_{X}) \cong
H^0(X,\sext^{r}_{\O_{X}}(\O_{Z},\Omega^{r}_{X}))$, we need only prove that
the classes agree locally.

Let $x \in X$ be a point and $U$ be a neighborhood of $x$ such that
the restriction of $E$ to $U$  is trivial, with $ \{ f_1, \ldots ,
f_r \} $ a local holomorphic framing of $E$ over $U$ and $ \{ f^1,
\ldots , f^r \}$ the dual framing for $E^{\vee}$; we identify $E$
and $E^{\vee}$ with the trivial bundle via these framings. Assume
that $Z$ has holomorphic equations $\{z_1, \ldots , z_r \}$ in $U$,
and is hence the zero set of the $\nu = z_1 f^1 + \cdots + z_r f^r$
of $E^{\vee}$. Then the Koszul complex $K(\nu)^{*}$ over $U$ is
quasi-isomorphic to $\O_{Z} {|}_{U}$ where $K(\nu)^{-i}=
\bigwedge^{i} \O_{U}^{r}$ and the differentials are contractions by
$\nu$ .

We will construct a quasi-isomorphism
$\tilde{u} : K(\nu)^* \to T^* |_{U}$ such that the class of the composition
$\Tr_{\Lambda}(\psi)|_{U} \circ \tilde{u}$ in
$$\sext^{r}_{\O_{X}}(K(\nu)^*,\A^{r,*}_{U}[r]) \cong
\sext^{r}_{\O_{X}}(\O_{Z}|_{U},\Omega^{r}_{U})$$ is the restriction of the
Grothendieck fundamental class of $Z$.

\textit{Step 1}:
We first define a map $u$ from $K(\nu)^{*}$ to $K(\tau)^{*}{|}_{U}$. By the
assumptions on $Z$ and $\tau$, there exist $u_{ji} \in \Gamma(U,\AX)$ such that
$z_i = \sum_j u_{ji}\alpha_j $ where $\tau= \alpha_1 e^1 + \cdots + \alpha_r
e^r$. We let $u_0 : \O_{U} \to  \A_U$ be the inclusion and $u_{-1}: \O^r_{U}
\to \A_{U} \otimes \E |_{U}$ be the map which sends $f_i$ to $\sum_j
u_{ji}e_j$. Then we extend $u_{-1}$ to a map of Koszul complexes by setting
$u_{-k} = \bigwedge^{k}u_{-1}$. Therefore $u_{-r}$ is given by multiplication by
$\det(u_{ij})$.

\textit{Step 2}: Next, we shall extend $u: K(\nu)^{*} \rightarrow
K(\tau)^{*}{|}_{U}$ to a map $ \tilde{u} : K(\nu)^{*} \rightarrow T^*{|}_{U}$.
The twisted complex $T^*$ is a filtered complex with respect to Dolbeault
degree. Let us denote this (decreasing) filtration by $F^i$, i.e. $F^i =
F^i(T^*)= \bigoplus_{j \geqslant i} \A^{0,j}_{X} \otimes \bigwedge \E$. Then one
has $Gr^i_F = \A^{0,i}_{X} \otimes K(\tau)^*$.

The map $(-\delta + \itau) : K(\tau)^* \to F^1 [1]$ is a map of complexes, since
$(-\delta + \itau) \circ \itau = -\delta \circ (-\delta + \itau)$. Hence
$$
(-\delta + \itau) \circ u : K(\nu)^* \to F^1 [1] |_{U}
$$
is a cochain map. Moreover, $(-\delta + \itau) \circ u$ is homotopic to the zero
map since $K(\nu)^*$ is a complex of free $\O_{X}$-modules; $F^1 [1]$ is acyclic
in negative degrees; and $\left((-\delta + \itau) \circ u \right) : \O_{U} \to
F^1 [1]^{0} |_{U}$ is the zero map. (Note that $u_0$ is the inclusion of
$\O_{U}$ into $\A_{U}$, and $(-\delta + \itau)_{0}=\dbar$). Consequently, there
exists an extension
$$
\tilde{u}=l(u): K(\nu)^* \to \cone(-\delta + \itau)^* [-1]
$$
by Lemma~\ref{la18}. This is the desired extension since $T^* =
\cone(-\delta + \itau)^* [-1]$.

\textit{Step 3}: We now prove that $\tilde{u} : K(\nu)^* \to T^* |_{U}$ is a
quasi-isomorphism. Let $i : \O_{Z} |_{U} \to \A^{0,*}_{U} \otimes \O_{Z}
|_{U}$ be the inclusion and $p: K(\nu)^* \to \O_{Z} |_{U}$ be the augmentation
map. Then we have a commutative diagram
$$
\CD
K(\nu)^* @>{\tilde{u}}>>   T^* |_{U} \\
@VpVV    @V{\epsilon}VV \\
\O_{Z} |_{U} @>i>>   \A^{0,*}_{U} \otimes \O_{Z} |_{U} \\
\endCD
$$
in which the vertical arrows and the bottom horizontal arrow are
quasi-isomorphisms. As a result, $\tilde{u}$ is a quasi-isomorphism.

\textit{Step 4}: Let $\eta = \Tr_{\Lambda}(\psi) \circ \tilde{u}$. Thus the degree
$(-r)$ component of $\eta$ is given by the composition $
\frac{1}{r!}(\intau)^{r} \circ \det(u_{ij})$. Hence
\begin{eqnarray*}
\eta_{-r} (f_1 \wedge \cdots \wedge f_r) &=& \frac{1}{r!}(\intau)^r(\det(u_{ij})e_1 \wedge
\cdots \wedge e_r ) \\
&=& \det(u_{ij})\frac{1}{r!}(\intau)^r(e_1 \wedge \cdots \wedge e_r) \\
&=& \det(u_{ij}) \partial \alpha_1 \wedge \cdots \wedge  \partial \alpha_r
\pmod{{\mathcal I}}  \\
&=& \partial z_1 \wedge \cdots \wedge \partial z_r \pmod{ {\mathcal I}} \\
&=& dz_1 \wedge \cdots \wedge dz_r \pmod{{\mathcal I}} \\
\end{eqnarray*}
By Proposition~\ref{la12} $\eta$ represents the Grothendieck
fundamental class of $Z \cap U$ in $U$. Since $\tilde{u}$ is a
quasi-isomorphism, so is $\Tr_{\Lambda}(\psi) |_{U}$. Since
$\Tr_{\Lambda}(\psi)$ represents the Grothendieck fundamental class
locally, it does so globally.
\end{proof}
\begin{lemma}   \label{la20}
Let $ \psi = \frac{1}{r!}R_{A}^{r}$. We have
$$
\Tr_{\Lambda}(\psi) {|}_{\AX} = \frac{1}{r!}\tr_{s}(\psi)
$$
\end{lemma}
\begin{proof}   Omitted.
\end{proof}
\begin{corollary}        \label{la21}
The image of the Grothendieck fundamental class in $H^r(X, \Omega^r_{X})$ is
represented by the $(r,r)$ degree part of the Chern character form of the
superbundle $\bigwedge E$ equipped with the superconnection $A = \nabla +
\delta$.
\end{corollary}
\begin{proof}
This follows from the theorem, Lemma~\ref{la20}, and Theorem~\ref{la6}.
\end{proof}
This completes the proof of Theorem B.
\nocite{Atiyah57}


\begin{thebibliography}{10}

\bibitem{Altman-Kleiman70}
{\scshape A.~Altman {\normalfont \smfandname} S.~Kleiman} -- \emph{Introduction
  to {G}rothendieck duality theory}, Lecture Notes in Mathematics, Vol. 146,
  Springer-Verlag, Berlin, 1970.

\bibitem{Atiyah57}
{\scshape M.~F. Atiyah} -- {\og Complex analytic connections in fibre
  bundles\fg}, \emph{Trans. Amer. Math. Soc.} \textbf{85} (1957), p.~181--207.

\bibitem{Atiyah-Hirzebruch62}
{\scshape M.~F. Atiyah {\normalfont \smfandname} F.~Hirzebruch} -- {\og
  Analytic cycles on complex manifolds\fg}, \emph{Topology} \textbf{1} (1962),
  p.~25--45.

\bibitem{Berline-Getzler-Vergne92}
{\scshape N.~Berline, E.~Getzler {\normalfont \smfandname} M.~Vergne} --
  \emph{Heat kernels and {D}irac operators}, Grundlehren der Mathematischen
  Wissenschaften, vol. 298, Springer-Verlag, Berlin, 1992.

\bibitem{Brown59}
{\scshape E.~H. Brown, Jr.} -- {\og Twisted tensor products. {I}\fg},
  \emph{Ann. of Math. (2)} \textbf{69} (1959), p.~223--246.

\bibitem{Dolbeault56}
{\scshape P.~Dolbeault} -- {\og Formes diff\'erentielles et cohomologie sur une
  vari\'et\'e analytique complexe. {I}\fg}, \emph{Ann. of Math. (2)}
  \textbf{64} (1956), p.~83--130.

\bibitem{Fischer76}
{\scshape G.~Fischer} -- \emph{Complex analytic geometry}, Springer-Verlag,
  Berlin, 1976, Lecture Notes in Mathematics, Vol. 538.

\bibitem{Griffiths-Harris78}
{\scshape P.~Griffiths {\normalfont \smfandname} J.~Harris} -- \emph{Principles
  of algebraic geometry}, Wiley Classics Library, John Wiley \& Sons Inc., New
  York, 1994, Reprint of the 1978 original.

\bibitem{Grothendieck57}
{\scshape A.~Grothendieck} -- {\og Th\'eor\`emes de dualit\'e pour les
  faisceaux alg\'ebriques coh\'erents\fg}, S\'eminaire Bourbaki, Vol.\ 4, Soc.
  Math. France, Paris, 1995, p.~Exp.\ No.\ 149, 169--193.

\bibitem{Kapranov99}
{\scshape M.~Kapranov} -- {\og Rozansky–-{W}itten invariants via {A}tiyah
  classes\fg}, \emph{Compositio} \textbf{115} (1999), p.~71--–113.

\bibitem{Quillen85}
{\scshape D.~Quillen} -- {\og Superconnections and the {C}hern character\fg},
  \emph{Topology} \textbf{24} (1985), no.~1, p.~89--95.

\bibitem{Toledo-Tong75}
{\scshape D.~Toledo {\normalfont \smfandname} Y.~L.~L. Tong} -- {\og The
  holomorphic {L}efschetz formula\fg}, \emph{Bull. Amer. Math. Soc.}
  \textbf{81} (1975), no.~6, p.~1133--1135.

\bibitem{Toledo-Tong76}
\bysame , {\og A parametrix for {$\overline \partial $} and {R}iemann-{R}och in
  \v {C}ech theory\fg}, \emph{Topology} \textbf{15} (1976), no.~4, p.~273--301.

\bibitem{Toledo-Tong78}
\bysame , {\og Duality and intersection theory in complex manifolds. {I}\fg},
  \emph{Math. Ann.} \textbf{237} (1978), no.~1, p.~41--77.

\end{thebibliography}

\providecommand{\bysame}{\leavevmode ---\ }
\providecommand{\og}{``}
\providecommand{\fg}{''}
\providecommand{\smfandname}{et}
\providecommand{\smfedsname}{\'eds.}
\providecommand{\smfedname}{\'ed.}
\providecommand{\smfmastersthesisname}{M\'emoire}
\providecommand{\smfphdthesisname}{Th\`ese}

\end{document}